\documentclass[5p,times]{elsarticle}
\RequirePackage{etex}
\usepackage{amssymb}
\usepackage{amsthm}
\usepackage{makecell}
\usepackage{amsmath}
\usepackage{color,xcolor}
\usepackage{algorithm, algpseudocode}
\newtheorem{remark}{Remark}
\usepackage[utf8]{inputenc}
\usepackage{booktabs, caption, makecell}

\usepackage{threeparttable}
\usepackage{bbm}
\usepackage{tikz}
\usetikzlibrary{positioning}

\usepackage{multirow}
\usepackage[figuresright]{rotating}
\usepackage{amsmath}
\newcommand{\R}{\mathbb{R}}
\usepackage{tablefootnote}
\newcommand{\Z}{\mathbb{Z}}
\DeclareMathOperator*{\minimize}{minimize}

\DeclareMathOperator*{\subt}{subject \ to}
\newtheorem{problem}{Problem}
\usepackage{tikz-network}

\usepackage{algorithm, algpseudocode}

\begin{document}

\begin{frontmatter}

\title{Taking the human out of decomposition-based optimization via artificial intelligence: Part I. Learning when to decompose}

\author{Ilias Mitrai}
\author{Prodromos Daoutidis\fnref{label2}} 
\fntext[label2]{Corresponding author}
\ead{daout001@umn.edu}
\address{Department of Chemical Engineering and Materials Science, University of Minnesota, Minneapolis, MN 55455}

\begin{abstract} 
In this paper, we propose a graph classification approach for automatically determining whether to use a monolithic or a decomposition-based solution method. In this approach, an optimization problem is represented as a graph that captures the structural and functional coupling among the variables and constraints of the problem via an appropriate set of features. Given this representation, a graph classifier is built to determine the best solution method for a given problem. The proposed approach is used to develop a classifier that determines whether a convex Mixed Integer Nonlinear Programming problem should be solved using branch and bound or the outer approximation algorithm. Finally, it is shown how the learned classifier can be incorporated into existing mixed integer optimization solvers. 
\end{abstract}

\end{frontmatter}

\section{Introduction}

Decision-making problems arise in a wide range of applications in chemical engineering, from the molecular to the enterprise scale \cite{hanselman2016mathematical, grossmann2012advances,daoutidis2018integrating, pistikopoulos2021process}. Mathematical programming has been widely used to model such decision-making (optimization) problems. Different solution methods have been developed for the solution of broad classes of optimization problems, such as Mixed Integer Linear (MILP) \cite{conforti2014integer}, Nonlinear (NLP) \cite{wachter2006implementation}, and Mixed Integer Nonlinear Programming (MINLP) problems \cite{tawarmalani2005polyhedral, misener2014antigone}. Despite these significant advances, the solution of optimization problems in chemical engineering applications is complicated by the nonlinear nature of most physical and chemical processes, the presence of multiple temporal and spatial scales, and the discrete nature of design and operational decisions. The combination of these features leads to complex and large-scale optimization problems whose monolithic solution is challenging. 
\begin{figure*}[h]
    \centering
    \includegraphics[scale = 0.5]{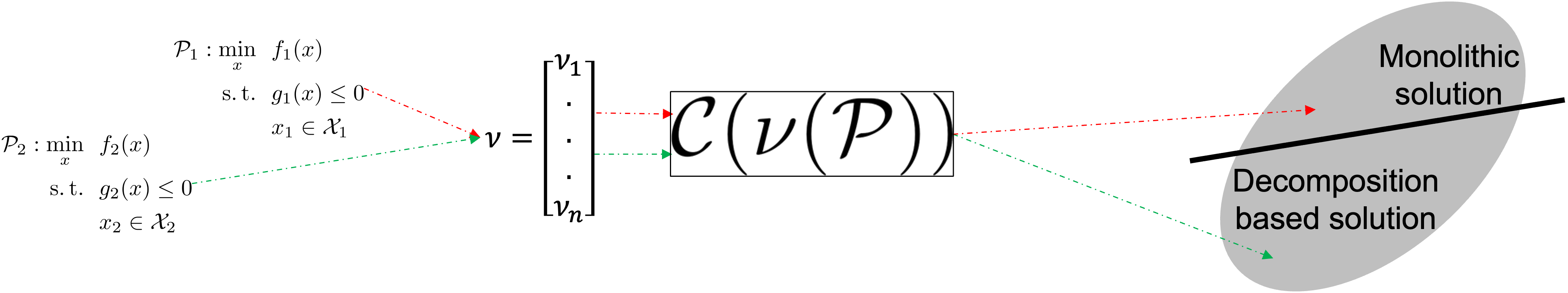}
    \caption{Algorithm selection via classification based on handcrafted features $\nu$}
    \label{fig:standard approach}
\end{figure*}
\begin{figure*}
    \centering
    \includegraphics[scale = 0.5]{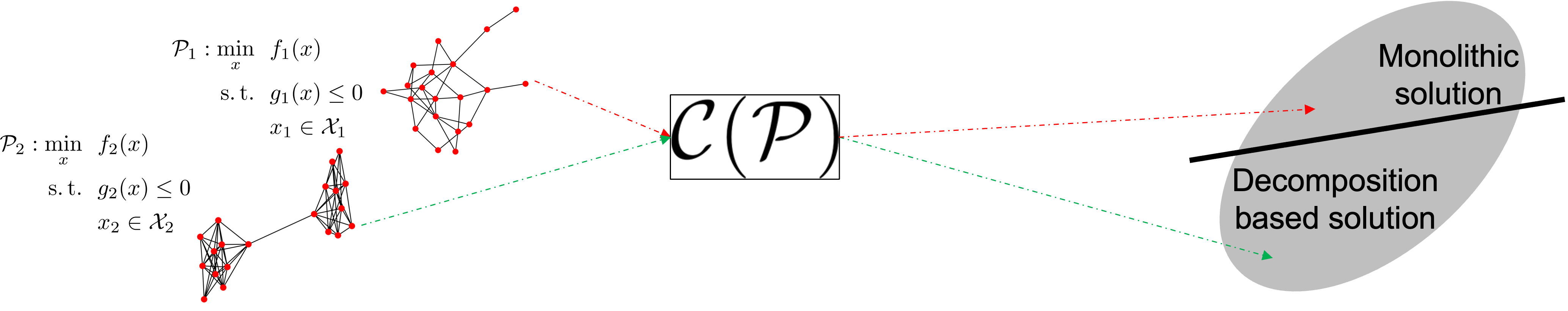}
    \caption{Algorithm selection via graph classification}
    \label{fig:proposed approach}
\end{figure*}

Decomposition-based solution algorithms have been widely used to solve such optimization problems by exploiting the underlying structure of the problem \cite{conejo2006decomposition,jiang2014optimal, li2009integrated,mouret2011new, aytug2015feature,bunel2020lagrangian,lara2018deterministic,tang2018optimal}. In general, these algorithms can be classified as distributed or hierarchical. In both cases, a complex optimization problem is decomposed into easier-to-solve subproblems. The main difference lies in the solution and coordination of the subproblems. In distributed algorithms, the subproblems are solved in parallel and are coordinated via dual variables. Typical examples include Lagrangean relaxation \cite{fisher2004lagrangian,fisher1985applications} and Lagrangean decomposition \cite{guignard1987lagrangean}, and the Alternating Direction methods of Multipliers \cite{boyd2011distributed}. In hierarchical algorithms, the subproblems are solved sequentially, based on some underlying hierarchy of the problem, and are coordinated via the addition of cuts. Typical examples are Benders \cite{benders1962partitioning}, Generalized Benders \cite{geoffrion1972generalized}, Outer Approximation \cite{duran1986outer, fletcher1994solving}, and Bilevel \cite{conejo2006decomposition} decomposition algorithms. Finally, there are algorithms that exploit both distributed and hierarchical structures such as cross-decomposition \cite{van1983cross}. 

Despite the advantages of decomposition-based solution methods their efficiency over monolithic solution approaches is not known a-priory since it depends on multiple factors. For example, the performance of Benders and Generalized Benders decomposition depends on the computational complexity of the master problem and the subproblem, the number of infeasible subproblems that are solved, and the quality of cuts that are generated. Similarly, for distributed algorithms such as Lagrangean decomposition, the convergence depends on the coordination scheme and the initialization of the dual variables. 

In general, state-of-the-art optimization algorithms are black-box systems and determining the computational performance of an algorithm for a given problem a-priory is difficult. In this spirit, the decision whether to use a decomposition or a monolithic-based solution algorithm for the solution of a given problem is not obvious in general. Although the number of variables and constraints can guide this determination, for some classes of problems such as convex MINLP problems decomposition-based methods can be more efficient than monolithic ones even for problems with few variables and constraints. Furthermore, even if such a determination is made, the implementation of the requisite method has multiple steps, such as problem decomposition, coordination, and initialization. The configuration of these steps is nontrivial and can have an important effect on the performance of the algorithm. Our overarching goal is to develop an automated framework to address these problems, effectively taking human judgment or intervention out of this process. In this paper, we focus on whether a problem should be solved using a monolithic or a decomposition-based solution strategy. In a companion (Part II) paper \cite{mitrai2023initviaAL}, we focus on the initialization of Benders decomposition for online applications.

The problem of selecting the best solution method for a given computational task is formally known as the algorithm selection problem \cite{rice1976algorithm, kerschke2019automated}. Although this is a well-studied problem in Computer Science it has received much less attention in the Process Systems Engineering community. The algorithm selection problem for optimization problems has three components: (1) the space of optimization problems $\mathcal{P}$, i.e., all the problems that can be considered, (2) the set of algorithms $\mathcal{A}$, i.e., all algorithms that are available for the solution of the problems in $\mathcal{P}$, and (3) the performance space $\mathcal{M}$, i.e., some metric used to compare solution methods for a given problem. Given these sets the algorithm selection problem is defined as follows:   
\begin{problem}
    \normalfont (Algorithm selection) Given an optimization problem $P$ and a set of algorithms $\mathcal{A} = \{a_1,...,a_n\}$, determine which algorithm $\alpha^*$ should be used to solve the problem such that a desired performance function $m: \mathcal{P} \times \mathcal{A} \rightarrow \mathcal{M}$ is optimized. 
\end{problem}
The performance function $m$ can be either the solution time or the best feasible solution for a given computational time budget. This problem can be posed as the optimization problem 
    \begin{equation} \label{alg sel}
        a^* \in \arg \min_{a \in \mathcal{A}} m(P,a).
    \end{equation}
The algorithm selection problem is a black-box optimization problem since the function $m$ is not known explicitly. Existing automated algorithm selection tools rely on Artificial Intelligence (AI) and Machine Learning (ML) to approximate the solution of the algorithm selection problem. A typical approach is to approximate the computational performance of an algorithm, measured in terms of solution time for global optimization solvers for a given problem (or class of problems), using a surrogate model \cite{smith2012measuring, hutter2014algorithm, kotthoff2016algorithm, gomes2001algorithm,pihera2014application}, such as Decision Trees \cite{bartz2004tuning}, Gaussian Processes \cite{hutter2006performance}, Neural Networks \cite{smith2011discovering}, and Ridge regression \cite{huang2010predicting}. The surrogate model is subsequently used to select the best solution method. Another approach is to approximate the solution of the algorithm selection problem itself using classification techniques, where a classifier is used to approximate the best solution strategy \cite{xu2011hydra}. These approaches have been applied to determine the best solution strategy for satisfiability problems, Mixed Integer Linear (MILP) and Mixed Integer Quadratic (MIQP) Programming problems  \cite{xu2011hydra,bonami2022classifier,kruber2017learning}. 

In these AI/ML approaches the optimization problem cannot be the direct input to the standard surrogate models mentioned earlier, since it is a complex data point having different types of variables and constraints and an objective function. Furthermore, two optimization problems can have different number and type of variables and constraints making their comparison difficult. This necessitates the usage of a handcrafted set of features (see Fig.~\ref{fig:standard approach}), such as the number of variables, constraints, range of coefficients in constraints and objective, the condition number of constraint matrix, etc. (see \cite{smith2012measuring,xu2011hydra,kruber2017learning} for a detailed description of such features). The usage of such features is essentially a dimensionality reduction step which although enables the application of standard ML tools, it also aggregates information about the problem and does not account for the coupling among the variables and constraints which is key in decomposition-based solution algorithms.

To overcome these limitations, in this paper we propose a new graph abstraction, suitable for generic nonlinear optimization problems, and a subsequent graph classification approach based on geometric deep learning to determine if a decomposition-based solution algorithm should be used. Specifically, based on our previous work we represent an optimization problem as a graph \cite{allman2019decode,mitrai2022stochastic}, where every node is a constraint or a variable. This representation captures the structural coupling among the variables and constraints of a problem. We extend this representation by adding features in every node. These features are obtained from the problem formulation and account explicitly for the functional form of the constraints; their concatenation forms the feature matrix. Thus, the structural coupling among the variables and the constraints is captured via the adjacency matrix and the functional coupling via the feature matrix. This graph representation allows us to adopt a geometric deep learning framework (graph neural networks) \cite{bronstein2021geometric} to build a graph classifier to approximate the solution of the algorithm selection problem (Eq.~\ref{alg sel}) as presented in Fig.~\ref{fig:proposed approach}. The proposed approach can be applied to generic optimization problems with different number and type of variables and constraints. For illustration, we apply the proposed approach to determine if a convex MINLP should be solved using Branch and Bound or the Outer Approximation algorithm. We use benchmark optimization problems for training and testing. The results show that the accuracy of the classifier on the testing dataset is 90$\%$. Finally, we show how the proposed approach can be integrated into existing software technology for the automated solution of convex MINLP problems. 

The rest of the paper is organized as follows: in Section \ref{graph repr} we present the proposed graph representation of an optimization problem; in Section~\ref{proposed framework} we pose the question of whether to decompose as an algorithm selection problem and we present the graph classification-based algorithm selection solution approach; in Section~\ref{application} we show how the proposed approach can be used to determine when to decompose convex MINLPs. 

\section{Graph representation of an optimization problem and node features} \label{graph repr}
Consider the following general optimization problem
\begin{equation} \label{full model}
	\begin{aligned}
			P := \minimize_{x,y} \ \ & f(x,y) \\
			\subt \ \ & g_i(x,y) \leq 0 \ \ \forall i=1,...,m_{in}\\
			& h_j(x,y) = 0 \ \ \forall j=1,...,m_{eq}\\
			& x \in \R^{n_x^c}, y \in \Z^{n_y^d},
		\end{aligned}   
\end{equation}
where $ n_x^c + n_y^d = n $, $ m_{in}+m_{eq}=m$. The computational complexity of this problem depends on the number of variables and constraints, the interaction pattern among the variables, and the form of the constraints and objective function, i.e., whether they are convex or nonconvex. 

\begin{figure*}[h]
    \centering
    \includegraphics[trim = 0 0 0 0, clip,scale = 0.4]{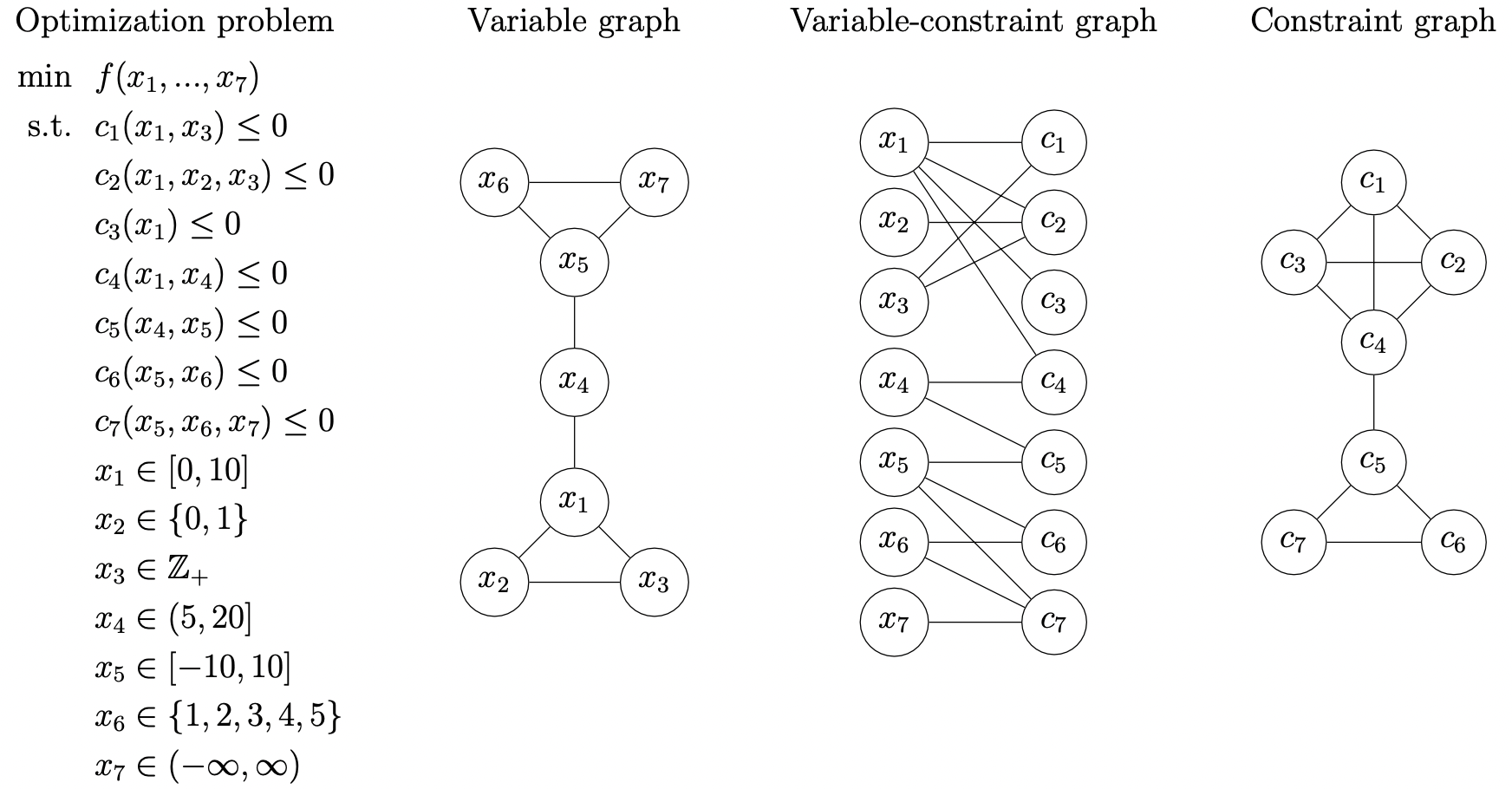}
    \caption{Graph representation of an optimization problem}
    \label{fig:graph representation}
\end{figure*}

In previous work, we have proposed a graph representation of optimization problems \cite{allman2019decode}. Three types of graphs were proposed (see Fig.~\ref{fig:graph representation}). The first is a bipartite variable-constraint graph, $\mathcal{G}_b(V_n, V_m, E)$ $(|V_n|=n, |V_m|=m)$, which has two sets of nodes, one representing the constraints $(V_m)$ and the other the variables $(V_n)$. The edges $E$ capture the presence of a variable in a constraint and form the adjacency matrix $A_b \in \mathbb{R}^{n \times m}$ with $A_{b}(i,j) \in \{0,1\}$, i.e., $A_b(i,j)=1$ if an edge exists between nodes $i$ and $j$, and zero otherwise. The second type of graph is the variable graph $\mathcal{G}_v(V_n, E_m)$ $(|V_n|=n)$ with adjacency matrix $A_v \in \mathbb{R}^{n\times n}$ where the nodes are the variables of the problem and the edges represent the constraints that couple the variables. The third type is the constraint graph $\mathcal{G}_c (V_m, E_m) (|V_m|=m, A_v \in \mathbb{R}^{m \times m})$ where the nodes are the constraints and the edges represent the variables that couple two constraints. This representation captures the structural coupling between the variables and the constraints, i.e. the presence or not of a variable in a constraint, as reflected in the adjacency matrix. Such graph representations have enabled the application of tools from network science and graph theory in order to ``learn'' the structure of an optimization problem which can be subsequently used as the basis for the application of decomposition-based solution algorithms \cite{allman2019decode, mitrai2022stochastic,mitrai2021iecr, mitrai2022multicut, mitrai2020decomposition}. Whereas a bipartite graph variable-constraint graph is the more general representation, the unipartite variable or constraint graphs may be better suited to specific decomposition-based solution methods depending on whether complicating variables or constraints are involved \cite{mitrai2022stochastic,mitrai2021iecr}. Yet, all these representations do not capture information about the type and domain of the variables or their functional form in the different constraints (linear, power law, etc.) in the case of nonlinear problems. We will refer to such information as functional coupling. 

\begin{figure*}
    \centering
    \includegraphics[trim = 0 0 0 0, clip,scale =0.25]{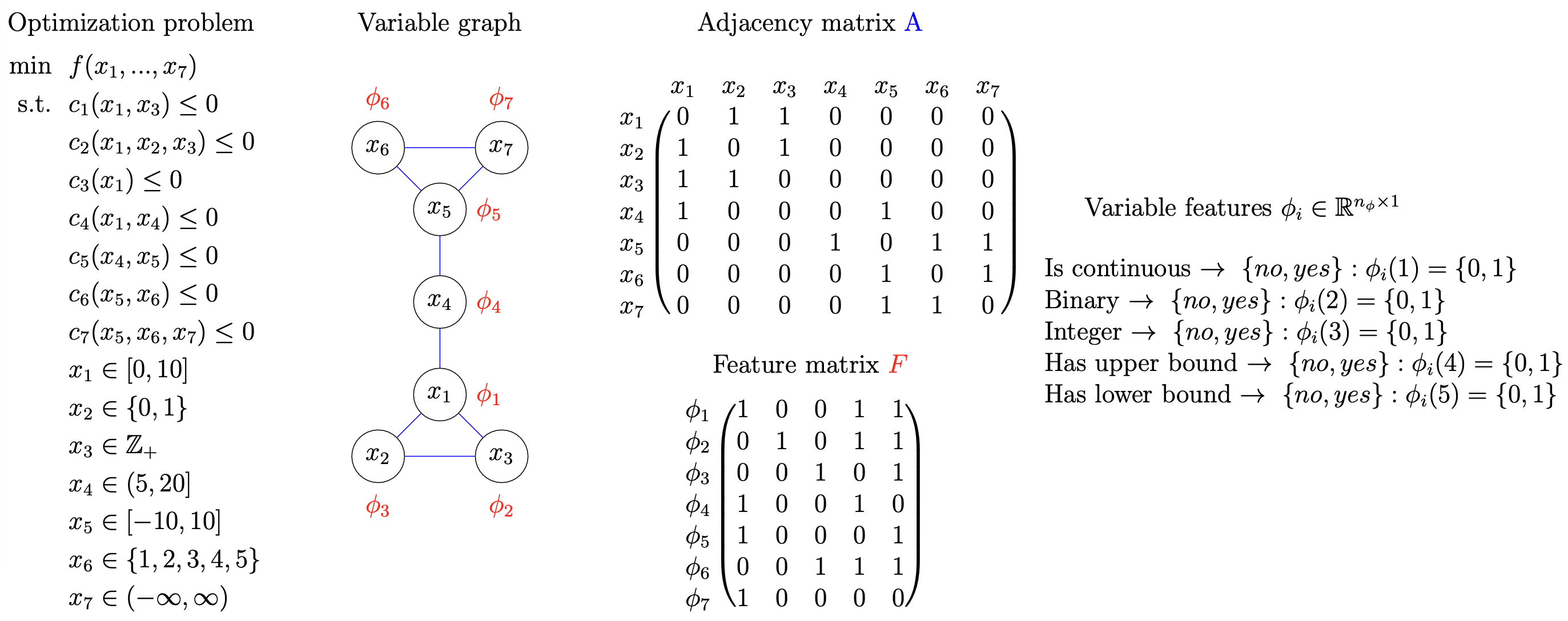}
    \caption{Graph and feature representation of an optimization problem}
    \label{fig:graph representation with features}
\end{figure*}

In this paper, we extend the aforementioned representation by incorporating a set of features for every node. There is considerable flexibility and generality in the type of features that can be added to capture different characteristics of the nodes. For example, such features can be the domain and the upper and lower bound for the variable graph, and the type of constraint (affine, convex, or nonconvex) for the constraint graph. Note that this representation can also capture information regarding the objective function of the problem since it can be considered as a constraint by reformulating the problem. For a given graph $G(V,E)$ with $|V|=n_V$ (this can be any graph presented in Fig.~\ref{fig:graph representation}), we can thus associate every node $n_i$ with a set of features $\phi_i \in \mathbb{R}^{1 \times n_{\phi}}$. Concatenation of these features forms the feature matrix $F \in \mathbb{R}^{n_V \times n_{\phi}}$. Under this representation, every graph is represented by three sets: the nodes $V$, the edges $E$, and the node features $F$, i.e., $G(V,E,F)$. This representation simultaneously captures the structural and functional coupling among the variables and the constraints of an optimization problem as presented in Fig.~\ref{fig:graph representation with features} (a variable graph is shown as an example). Finally, one can also consider features for the edges that couple different variables/constraints. For example, for the bipartite graph the edge features can represent the functional form of a variable in a constraint; in the variable graph, the edge features can represent the functional form in which the variables are present in the constraints that couple them; and in the constraint graph the edge features can represent the type of variables that couple two constraints. These features can in turn form an edge features matrix. We note that similar ideas have been proposed in the literature \cite{ding2020accelerating,gupta2020hybrid,gupta2022lookback,li2022learning, liu2022learning,nair2020solving,paulus2022learning}, however, they have been applied to MILP problems and for tuning rather than selecting optimization solvers.

\section{Learning when to decompose via graph classification} \label{proposed framework}

In this section, we present the graph classification framework which will approximate the solution of the algorithm selection problem presented in Eq.~\ref{when to dec alg sel}. For illustration, we use the variable graph of an optimization problem. 

\subsection{Learning when to decompose as an algorithm selection problem}\label{learn when to dec as alg sel prob}
First, we pose the question of whether to decompose as an algorithm selection problem. Specifically, the set $\mathcal{P}$ represents the class of optimization problems which are considered, i.e., MILP, convex MINLP, etc. The set $\mathcal{A}$ has two algorithms, a monolithic-based one $\alpha_M$ and a decomposition-based one $\alpha_D$, $\mathcal{A} = \{\alpha_M, \alpha_D\}$. Finally, the performance space is the solution time, i.e., $\mathcal{M} = \mathbb{R}_{+}$. Given these sets, we can determine whether to solve a given problem using a decomposition-based or a monolithic solution algorithm by solving the following problem:
    \begin{equation} \label{when to dec alg sel}
        a^* \in \arg \min_{a \in \{\alpha_M, \alpha_D\}} m(P,a).
    \end{equation}

In general, multiple monolithic and decomposition-based solution algorithms might be available. This scenario can be easily accommodated by increasing the size of the set of possible algorithms $\mathcal{A}$. 
Also, in the above formulation the performance function was chosen as the solution time. Depending on the problem of interest, one can select a different metric, such as the duality gap or the quality of the best feasible solution found after a fixed computational budget.
\begin{remark}
\normalfont The question of whether to use a decomposition-based solution algorithm over a monolithic one can be posed for every class of problems. However, different performance functions might be appropriate for different classes of problems and algorithms, based on the available convergence guarantees. For example, one can consider whether a nonconvex MINLP should be solved using Branch and Bound or Generalized Benders Decomposition (GBD). Since the convergence of GBD is not guaranteed for this class of problems, the performance function cannot be simply the solution time, since GBD can converge faster than branch and bound but the solution can be highly suboptimal.
\end{remark}

\subsection{Graph classification approach and architecture}
\begin{figure*}
    \centering
    \includegraphics[trim = 0 0 0 0, clip,scale=0.37]{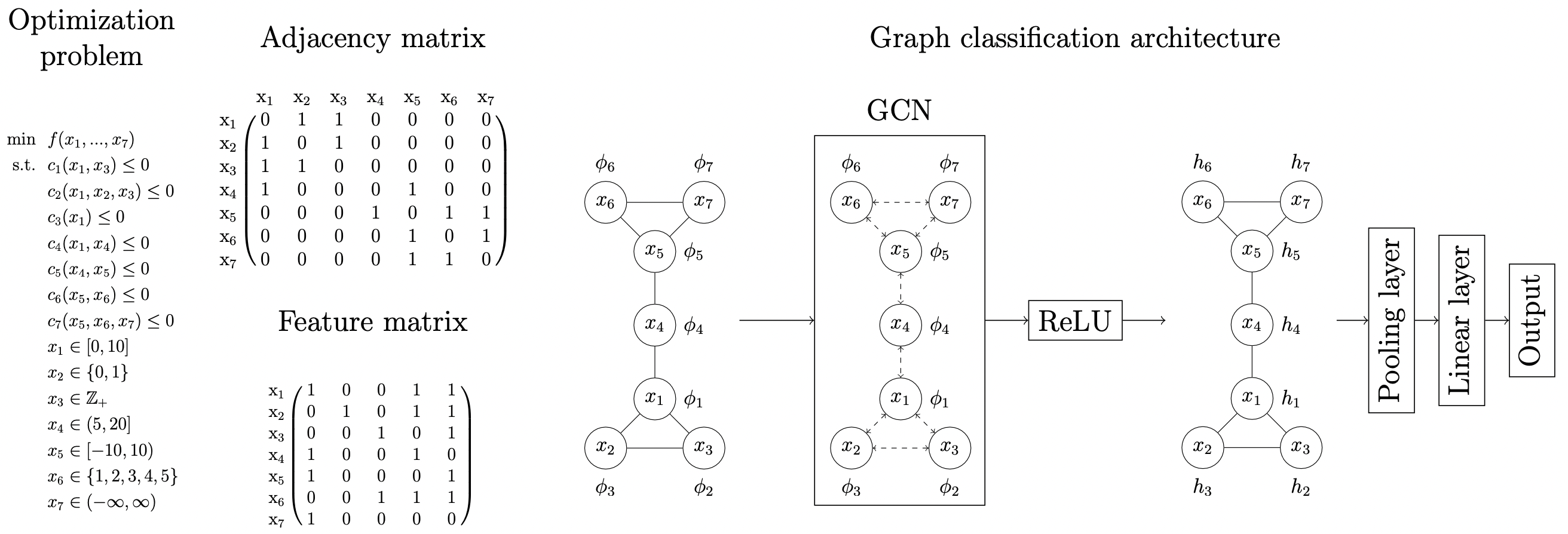}
    \caption{Learning when to decompose framework}
    \label{fig:framework}
\end{figure*}

Given an optimization problem P and the variable graph $\mathcal{G}_n (V, F, E)$ with adjacency matrix $A$, the goal is to develop a classifier $\mathcal{C}: F(\mathcal{P}) \times A(\mathcal{P}) \mapsto p \in \mathbb{R}^{N_a}$ to determine if problem P should be solved using a monolithic $(\alpha_M)$ or decomposition-based $(\alpha_D)$ solution approach ($\mathcal{A}=\{\alpha_M, \alpha_D\},|\mathcal{A}|=N_a=2)$. 
The inputs in the classifier are the adjacency $A(P)$ and feature $F(P)$ matrices (which depend on the problem $P$), and the output is a vector $p \in \mathbb{R}^{N_a \times 1}$, where $p_i$ is equal to the probability that algorithm $\alpha_i$ solves the problem in the minimum computational time. 
Under this setting, the algorithm selection problem is transformed into a graph classification problem, where for a given problem $P$ the best solution strategy $\alpha^*$ is 
\begin{equation} \label{learn when to dec prop eq}
    \alpha^* = \arg \max \ \{p_i\}_{i=1}^{N_a}, 
\end{equation}
where with $p=\mathcal{C}(F(P), A(P))$. Comparing Eq.~\ref{when to dec alg sel} with Eq.~\ref{learn when to dec prop eq} we observe that the graph classifier $\mathcal{C}(F(P), A(P))$ approximates $ \min_{a \in \{\alpha_M, \alpha_D\}} m(P,a)$ and selects the algorithm which has the highest probability to solve problem $P$ in the minimum computational time. The overall framework is presented in Fig.~\ref{fig:framework}. 

The prediction of the best solution strategy is performed by considering the exact structural and functional coupling among the variables and the constraints of the problem. This prediction has three steps: (1) message passing, (2) pooling, (3) final classification step. The first step updates the features of a node by considering the features of the adjacent nodes, propagating information about the features of a node across the graph. The pooling layer, creates a graph-level feature which is used to characterize the whole graph based on the features of the individual nodes. The last step performs the classification step based on the graph level features. In the rest of this section we present each step in detail using the variable graph as an example. However, both the constraint and the variable-constraint bipartite graphs can be used. 

\subsection{Message passing}
Given the graph representation of the problem and the adjacency $A$ and feature matrices $F$, message passing is performed to update the features of the nodes using information from the neighbors. Specifically, given an optimization problem $P_i$ (see Eq. 2) with $N_v$ variables and $N_m$ constraints, we can generate the variable graph $G_v  (V_v,E_v) (|V_v |=N_v)$ and obtain the adjacency $A \in \mathbb{R}^{N_v \times N_v}$ and feature $F \in \mathbb{R}^{N_v \times N_\phi}$ matrices where $N_\phi$ is the number of features per node. We define as $\phi_i$ the features of node i ($\phi_i$ is the $i^th$ row of matrix F). The features of node $i$ are updated using information about the features of the neighbor nodes $\mathcal{N}(i)$ of node $i$ as presented below \cite{bronstein2021geometric}
\begin{equation}
    \phi_i^{updated} = \sigma \bigg( \phi_i , \bigoplus_{j \in \mathcal{N}(i)} \psi(\phi_i, \phi_j;W) \bigg), 
\end{equation}
where $\phi_i^{updated}$ are the updated features of node $i$, $\bigoplus$ denotes an aggregation function that is independent of the order of the neighbors (sum, average, max, etc.), $\psi$ is a message function that takes as input the features of node $i$ and the neighbor $j\in \mathcal{N}(i)$, $W$ as learnable weights, and $\sigma$ is an update function that returns the new features of node $i$. Based on the type of function $\psi$ and the aggregation function $\bigoplus$ that are used different graph neural network models have been proposed, such as Graph Convolutional Neural Networks (GCN) \cite{kipf2016semi}, Graph Attention Networks (GAT) \cite{velivckovic2017graph} and GraphSage \cite{hamilton2017inductive} (see \cite{bronstein2021geometric} for a review).

\begin{figure}
    \centering
    \includegraphics[trim = 0 0 0 0, clip,width=\columnwidth]{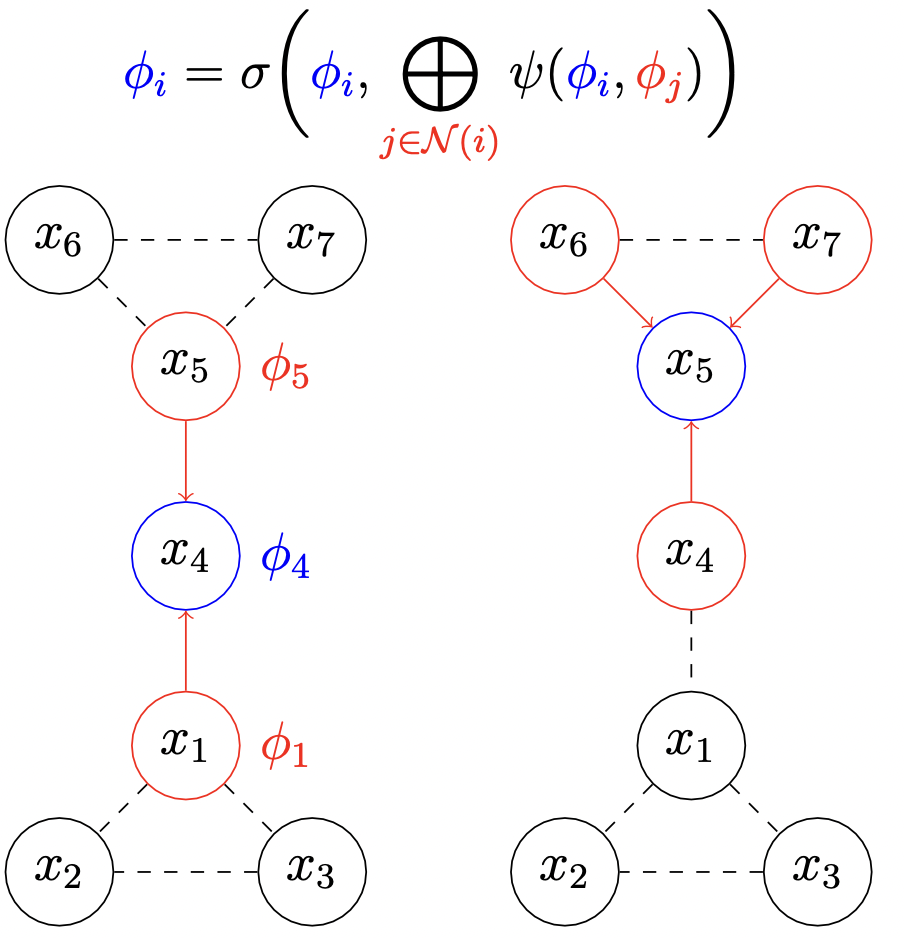}
    \caption{Message passing on the graph representation of the optimization problem}
    \label{fig:message passing opt prob}
\end{figure}

Stacking $L$ such layers sequentially leads to a deep learning architecture where the features of a node in layer $l$ are given by:
\begin{equation}
\begin{aligned}
    h^{1}_i   & = \phi_i \\ 
    h_i^{l+1} &= \sigma \bigg( h_i^{l} , \bigoplus_{j \in \mathcal{N}(i)} \psi(h_i, h_j;W^l) \bigg) \ \forall l=1,..,L-1,  
\end{aligned}
\end{equation}
where $h_i^l \in \mathbb{R}^{N_h \times 1} \ \forall l\geq 1$ are the features of node $i$ in level $l$, $W^l$ are the learnable weights at layer $l$, and $N_h$ is the dimension of the hidden features. An example of the message passing on the variable graph of the optimization problem is shown in Fig~\ref{fig:message passing opt prob}.

\subsection{Pooling}
The output of the last message passing layer is the original graph but the features of the nodes have been updated to $H$. At this stage, the feature matrix $H$ can not be used for classification, since the dimensions of $H$ depend on the number of the nodes (variables and constraints) in the graph of the optimization problem. To overcome this, a pooling layer is used to create a graph-level feature $r \in \mathbb{R}^{N_{\phi}}$. Different pooling functions can be used \cite{bronstein2021geometric,grattarola2022understanding}, however, in this paper $r$ is equal to
\begin{equation}
r_i   =\frac{1}{N_v}  \sum_{j=1}^{N_v} H_{ij}^L   \ \  \forall i=1,...,N_v.    
\end{equation}

This pooling function returns the average value of every feature across all the nodes in the graph. 
\subsection{Final classification step}
Finally, the graph level feature is an input to a linear transformation layer, where the output is $ y \in \mathbb{R}^{N_a \times 1}$ ($y_i$ is the probability that algorithm $i$ solves the problem in the minimum computational time) and is equal to 
\begin{equation}
y=\Theta r+b,    
\end{equation}
where $\Theta \in \mathbb{R}^{N_A \times N_h}$, $b \in \mathbb{R}^{N_A \times 1}$ are parameters. 

\begin{table*}[ht]
\centering
\caption{Performance metrics for the graph classifier on the testing data set}
\begin{tabular}{ccccccc}
\hline
\multirow{2}{*}{Accuracy} & \multirow{2}{*}{0.9} & \multirow{2}{*}{} & \multirow{2}{*}{} & \multicolumn{3}{c}{Confusion matrix}                  \\ \cline{5-7} 
                          &                      &                           &                      &                 & \multicolumn{2}{c}{Predicted label} \\ \hline
Label                     & Precision            & Recall                    & F1 score             & True label      & OA               & B\&B             \\ \hline
OA                        & 0.83                 & 1.00                      & 0.91                 & OA              & 15               & 0                \\
B\&B                      & 1.00                 & 0.80                      & 0.89                 & B\&B            & 3                & 12               \\ \hline
\end{tabular} 
\label{table: classifier results}
\end{table*}
\subsection{Training}
We assume that a set of optimization problems $\{P_i\}_{i=1}^{N_{data}}$ is available and for these problems we can obtain the variable graphs $\{G_i\}_{i=1}^{N_{data}}$, and the adjacency and feature matrix $\{A_i,F_i\}_{i=1}^{N_{data}}$ and labels $\{z_i\}_{i=1}^{N_{data}}$ where $z_i$ is the algorithm that solves problem $P_i$ in the minimum computational time. Given these data, the task is to learn the classifier's parameters that will maximize its accuracy by optimizing some loss function, such as cross-entropy. The learning problem is
\begin{equation}
    \begin{aligned}
        \minimize_{ \{W_i\}_{i=1}^{L-1}, \theta, b} \ \ \mathcal{L} \bigg( z^{data}, z^{pred};\{W_i\}_{i=1}^{L-1},\theta,b \bigg)
    \end{aligned}
\end{equation}

\begin{remark}
\normalfont This classifier can classify problems with different number of variables and constraints since the pooling layer after the last convolution layer creates the vector $r$ which has $N_h \times 1$ dimension. Therefore, for the final classification step, every problem is represented in the hidden features dimension space $N_h$.
\end{remark}
\begin{remark}
\normalfont 

In the case where multiple monolithic and decomposition solution methods are available, the same architecture can be used, however, the resulting problem becomes a multiclass classification one. 
\end{remark}

\section{Application to convex MINLP problems} \label{application}
In this section we implement the proposed approach to determine whether a convex MINLP problem should be solved using Branch and Bound (B\&B) \cite{gupta1985branch}, a monolithic-based solution approach, or the Outer Approximation (OA) algorithm \cite{duran1986outer}, a decomposition-based solution approach, as implemented in BONMIN \cite{bonami2008algorithmic}. 

\subsection{Branch and Bound}
Branch and Bound can be used to solve a convex MINLP problem by branching on the binary variables. Given a convex MINLP as presented in Eq.~\ref{full model}, branch and bound starts by solving the continuous relaxation of the problem which is a convex NLP. We will denote the solution of this problem as $x^{NLP},y^{NLP}$. Given this solution, branching on the binary variables is performed using some branching rules to select a binary variable $y_i$ and create two nodes in the branch and bound tree, one solved for $y_i$ fixed to $\lfloor y_i^{NLP} \rfloor$ and the other for $y_i$ equal to $\lceil y_i^{NLP} \rceil$. This procedure continues until the global optimal solution of the problem is found. 
 
\subsection{Outer Approximation algorithm}
The Outer Approximation algorithm \cite{duran1986outer} is a decomposition-based solution algorithm that alternates between the solution of a MIP master problem and a nonlinear convex optimization problem. Given a value of the binary variables, we obtain the subproblem
\begin{equation}
    \begin{aligned}
        \minimize_{x} \ \ & f(x,\bar{y}) \\
        \subt \ \ & g(x,\bar{y}) \leq 0\\
        & x \in \mathbb{R}^{n_x^c}.
    \end{aligned}
\end{equation}
Given the optimal solution of the subproblem $\bar{x}$ for a fixed value of $y=\bar{y}$, the constraints and objective of the subproblem are approximated via hyperplanes as follows
\begin{equation}
    \begin{aligned}
        & f(x,y) \geq f(\bar{x},\bar{y}) + \nabla f(\bar{x},\bar{y})^\top \begin{bmatrix} x-\bar{x} \\ y-\bar{y} \end{bmatrix}  \ \forall x \in \mathbb{R}^{n}, y \in \mathbb{Z}^{n_y}\\
        & g(\bar{x},\bar{y}) + \nabla g(\bar{x},\bar{y})^\top \begin{bmatrix} x-\bar{x} \\ y-\bar{y} \end{bmatrix} \leq 0 \ \ \forall x \in \mathbb{R}^{n_x^c}, y \in \mathbb{Z}^{n_y^d}.
    \end{aligned}
\end{equation}
These approximations are generated iteratively as dictated by the master problem which is equal to
\begin{equation}
    \begin{aligned}
        \minimize_{x,y} \ \ & \eta  \\
        \subt \ \ & \eta \geq f(\bar{x}^q,\bar{y}^q) + \nabla f(\bar{x}^q,\bar{y}^q)^\top \begin{bmatrix} x-\bar{x}^q \\ y-\bar{y}^q \end{bmatrix} \ \forall q \in \mathcal{Q} \\
        & g(\bar{x}^q,\bar{y}^q) + \nabla g(\bar{x}^q,\bar{y}^q)^\top \begin{bmatrix} x-\bar{x}^q \\ y-\bar{y}^q \end{bmatrix} \leq 0 \ \forall l \in \mathcal{Q}\\
        & x \in \mathbb{R}^{n}, y \in \mathbb{Z}^{n_y},
    \end{aligned}
\end{equation}
where $\mathcal{Q} = \{1,..,N_{\mathcal{Q}}\}$ denotes iteration number. The OA algorithm alternates between the solution of the master problem, which provides a lower bound and the subproblem which provides an upper bound. This is the standard application of the algorithm assuming the subproblem is always feasible (see \cite{duran1986outer} for more details).

\subsection{Feature representation of the problem}
We will consider the variable graph of the optimization problem and the following features per node:
\begin{itemize}
    \item Variable domain: continuous, binary, integer 
    \item Upper bound
    \item Lower bound
\end{itemize}
We use one hot encoding to represent these features as follows:
\begin{itemize}
    \item $\phi_1 \in \{0,1\}$: 1 if the variable is continuous and 0 otherwise
    \item $\phi_2 \in \{0,1\}$: 1 if the variable is binary and 0 otherwise
    \item $\phi_3 \in \{0,1\}$: 1 if the variable is integer and 0 otherwise
    \item $\phi_4 \in \{0,1\}$: 1 if the variable has an upper bound and 0 otherwise
    \item $\phi_5 \in \{0,1\}$: 1 if the variable has a lower bound and 0 otherwise
\end{itemize}
Using these features the dimensions of the feature matrix are $N_v \times 5$. 

\subsection{Data gathering for classification}

We use benchmark convex MINLP problems to train and test the classifier. For every problem $P_i$, we obtain the adjacency matrix $A_i$ using DecODe \cite{mitrai2022stochastic} and the feature matrix $F_i$ using the procedure presented in Algorithm~\ref{alg:get feat of var}. Every problem is solved with both algorithms  using BONMIN 1.8.8 \cite{bonami2008algorithmic} with a maximum computational time of 3000 seconds. All the other parameters are set equal to their default values. From the 295 problems, 227 are solved with at least one solver (151 problems are solved faster with OA and 57 with B\&B). Given the solution times, we use Algorithm~\ref{alg:get labels} to obtain the label $z_i$ of every problem. Overall, we obtain the dataset $\mathcal{D}=\{(A_i,F_i),z_i \}_{i=1}^{N_{data}}$, where $66\%$ of the data points have label OA and $34\%$ B\&B. 

\subsection{Graph classification architecture and implementation}
We split the dataset $\mathcal{D}$ at random into a training set and a testing set. The training set has 197 data points and the testing set has 30 data points picked at random (15 random data points have label OA and 15 have label B\&B). We perform hyperparameter optimization regarding the number of hidden features ($N_{\phi} = \{12$,$16$,$24$,$32$,$64$,$128\}$), the batch size ($\{10,20,40,50\}$), and the learning rate ($\{0.0001$, $0.0005$, $0.001$, $0.005$, $0.01$, $0.05\}$).

Based on the hyperparameter optimization results, the graph classifier has four convolution layers $(L=4)$, a (global mean) pooling layer, a linear layer (with dropout probability equal to 0.5), the number of hidden features is 12 $(N_h=12)$, and \verb|tahnh| is used as the activation function. For the training, we use the Adam algorithm \cite{kingma2014adam} for 50 iterations with random initialization, the learning rate equal to 0.005, and batch size equal to 10. The loss function is the cross entropy and different weights are assigned to the two classes to account for the imbalance in the training dataset. The weights are computed as follows: 
\begin{equation}
\omega_i= \frac{N_{data}}{N_\alpha} \sum_{j=1}^{N_{data}} \mathbbm{1}_i (z_j ))^{-1},    
\end{equation}
where $\mathbbm{1}_i$ is the indicator function, i.e., $\mathbbm{1}_i (z_i) =1$ if $z_j=i$ and 0 otherwise. The weight for class B\&B is 1.6148 and for class OA is 0.7243. The GCN is implemented in PyTorch Geometric \cite{fey2019fast} and the training is done using PyTorch \cite{paszke2019pytorch}.

\subsection{Graph classification results}
The accuracy of the classifier on the testing dataset is presented in Table~\ref{table: classifier results}. From the results we observe that the accuracy of the classifier is $90\%$ and all the problems with label $OA$, i.e., solved faster with the outer approximation algorithm, are classifier correctly and only 3 problems solved faster with branch and bound as mis-classified. These results can be attributed to the small dataset since only 197 data points are used for learning and the imbalance in the training set. 
\begin{figure*}
    \centering
    \includegraphics[trim = 10 0 0 0, clip,scale=0.4]{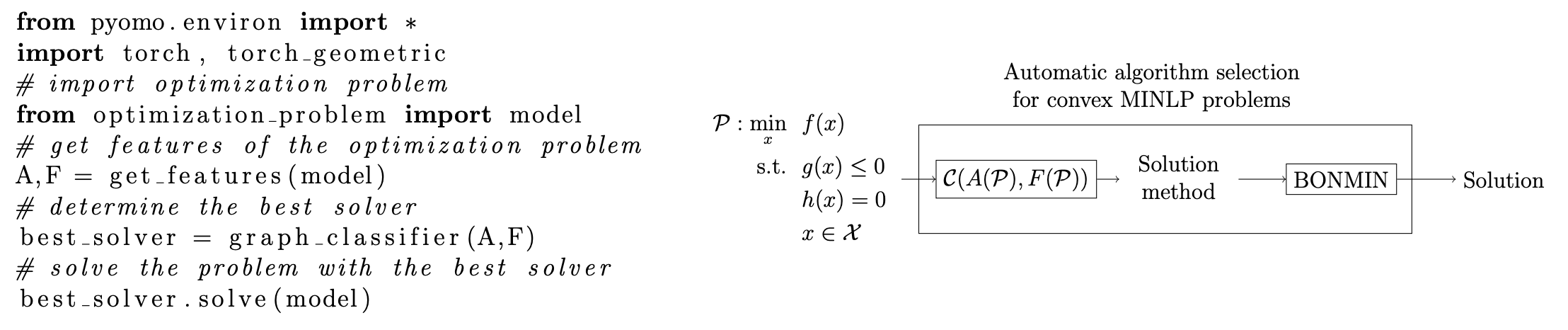}
    \caption{Automated algorithm selection for convex MINLP problems}
    \label{fig: aut alg sel convex minlp}
\end{figure*}
\subsection{Automated algorithm selection for convex MINLP problems}
The classification results show that the proposed approach was effective for a convex MINLP problem. The graph classifier can be readily incorporated in current mixed integer optimization software technology. Specifically, we consider the integration of the graph classifier with BONMIN and Pyomo \cite{hart2017pyomo} as presented in Fig.~\ref{fig: aut alg sel convex minlp}. Regarding the implementation in Pyomo, the graph classifier will predict the best solution strategy $\alpha^*$ which is either \verb|B-BB| (Branch and Bound) or \verb|B-OA| (Outer Approximation). The best solution strategy is passed to BONMIN as  \verb|solver.options['algorithm']=|$\alpha^*$.

\section{Conclusions and discussion}
Decomposition-based solution algorithms have been widely used to solve complex and large optimization problems. However, their efficiency over monolithic methods is not known a-priory. In this paper, we approached this problem as an algorithm selection one, and we proposed a new graph abstraction of an optimization problem and a subsequent graph classification approach for its solution. The proposed graph abstraction can account for the detailed structural and functional coupling among the variables and constraints of the problem, alleviating the need for a set of handcrafted features employed in previous approaches. The proposed approach was applied successfully to the problem of determining whether Branch and Bound or the Outer Approximation method should be used for the solution of convex MINLP problems. The proposed graph classifier can be easily integrated with optimization solvers, leading to an automated algorithm selection framework for decomposition-based optimization. 

The proposed approach requires data regarding the solution time or solution quality of optimization problems with different solvers. Such data can be obtained using the numerous available benchmark optimization problems \cite{gleixner2021miplib,minlplibANDnlp} and continuously enriched with data obtained from new problems. The availability of such data will enable the application of this approach to other classes of problems such as LP, MILP, and convex MIQPs in the future.

\section{Acknowledgements}
Financial support from NSF-CBET is gratefully acknowledged. 

\begin{algorithm}[h]
\caption{Procedure to obtain the label of an optimization problem}
\label{alg:get labels}
\begin{algorithmic}[1]
\Require{Optimization problem $\mathcal{P}$, Outer approximation algorithm (\verb|OA|), Branch and Bound algorithm (\verb|BB|), time limit (\texttt{time\_limit})} 
\For{$i \in \{\ OA,BB\}$ }
\State{Solve the optimization problem using algorithm $i$ with time limit equal to \texttt{time\_limit}}
\State{Store the solution time $t_i$}
\EndFor
\State{Find solver with minimum CPU time $i^* = \arg \min_{i \in \{OA,BB\}} t_i$}
\If{$i^* <$ \texttt{time\_limit}}
\If{$i^*=\texttt{BB}$}
\State{$y=1$}
\EndIf
\If{$i^*=\texttt{OA}$}
\State{$y=0$}
\EndIf
\Else
\State{\Return Problem not solved, data point not considered}
\EndIf
    \State \Return {Label $y$}
\end{algorithmic}
\end{algorithm}

\begin{algorithm}[h]
\caption{Procedure to obtain the features of a variable}
\label{alg:get feat of var}
\begin{algorithmic}[1]
\Require{Variable $v$} 
\State{Initialize the feature vector $\phi=[ 0\ 0\ 0\ 0\ 0]$}
\If{ $v$ is continuous}
\State{Set $\phi(1)=1$}
\State{Set $\phi(2)=0, \phi(3)=0$}
\Else
\State{Set $\phi(1)=0$}
\EndIf
\If{ $v$ is binary}
\State{Set $\phi(2)=1$}
\State{Set $\phi(1)=0, \phi(3)=0, \phi(4)=1, \phi(5)=1$}
\Else
\State{Set $\phi(2)=0$}
\EndIf
\If{ $v$ is integer}
\State{Set $\phi(3)=1$}
\State{Set $\phi(1)=0, \phi(2)=0, \phi(4)=1, \phi(5)=1$}
\Else
\State{Set $\phi(2)=0$}
\EndIf
\If{$v$ has upper bound}
\State{Set $\phi(4)=1$}
\Else
\State{Set $\phi(4)=0$}
\EndIf
\If{$v$ has lower bound}
\State{Set $\phi(5)=1$}
\Else
\State{Set $\phi(5)=0$}
\EndIf
\State{Get features $\phi = [\phi(1) \ \phi(2) \ \phi(3) \ \phi(4) \ \phi(5)] $} 
    \State \Return {Features $\phi$}
\end{algorithmic}
\end{algorithm}

\bibliographystyle{elsarticle-num}
\bibliography{sample}

\end{document}